# Fleet-mix Electric Vehicle Routing Problem for the E-commerce Delivery with Limited Off-Hour Delivery Implementation

Hyun-Seop Uhm[1,*], Abdelrahman Ismael[1], Natalia Zuniga-Garcia[1], Olcay Sahin[1], James Cook[1], Joshua Auld[1], Monique Stinson[2]

**ABSTRACT**

Freight truck electrification for last-mile delivery is one of the most important research topics to reduce the dependency on fossil fuel operations. Although a battery electric truck still has limitations on daily operations with lower driving ranges and higher purchasing cost than a conventional truck, operations with electrified trucks reduce total energy usage and driving noise on routes. In this paper, we propose a fleet-mix and multi-shift electric vehicle routing problem for joint implementation of fleet electrification and off-hour delivery in urban e-commerce delivery systems. Every electrified truck is assumed to have two shifts for both daytime and nighttime delivery operations while conventional trucks can operate during daytime only because of municipal restrictions on nighttime deliveries, which are related to engine noise. Also, every electrified truck must recharge between shifts at its depot. A fleet owner decides the best electrification ratio of the fleet and the proper number of chargers which gives the minimum total cost. The optimization problem is described as a mixed-integer linear programming model including common constraints for vehicle routing problem, recharging constraints, and two-shift operation of electrified trucks. A bi-level VNS-TS heuristic is also suggested for efficient solution search. The upper-level problem assigns trucks with engine type and brief route information using variable neighborhood search heuristic, and the lower-level problem finds the best route of each assigned truck using a tabu search heuristic. Scenarios with different EV driving ranges and nighttime operation availabilities are developed and evaluated with the POLARIS transportation simulation framework, and results are reported.

---

[1] TPS Transportation Systems & Mobility, Argonne National Laboratory, Lemont, IL, 60439
* Corresponding author: huhm@anl.gov

[2] U.S. Department of Transportation, Bureau of Transportation Statistics, Washington, DC, 20590
   Work conducted while Dr. Stinson was at Argonne National Laboratory

# 1. INTRODUCTION

The widespread use of fossil fuels has caused greenhouse gas (GHG) emissions, which have a negative impact on climate and environment *(1)*. The transportation sector is the largest contributor to GHG emissions in the United States, accounting for 27% of the total GHG emissions in 2020, of which 26% are from heavy-duty (HD) and medium-duty (MD) trucks *(2)*. Several research studies have proposed alternative methods to reduce truck emissions, such as: mechanical improvement of internal combustion engine (ICE) trucks, modal shift from trucks to more energy-efficient modes (e.g., rail), and truck electrification *(3)*. Promoting battery electric trucks can reduce GHG emissions significantly when the electricity is generated from a variety of clean and renewable sources *(4,5)*. Many conventional trucks that operate daily more than 100 miles, 10 hours, and 120 stops in urban areas for e-commerce deliveries could also reduce overall energy usage if they are electrified.

Fleet owners are facing pressure to electrify their trucks due to regulations on carbon-emission and the energy efficiency of electrified trucks, however, truck electrification currently faces three implementation barriers, including (a) low driving range, (b) high purchasing price of electrified trucks and electric vehicle supply equipment (EVSE), and (c) lack of public recharging infrastructure. Electric truck manufacturers aim for a driving range of 200–240 miles for their MD electric models, which is about 43% of the average driving range of current conventional trucks. A number of EV charging plugs for passenger cars have been recently installed at malls, marts, and other public parking spaces, however, these EVSEs are not suitable for truck recharging due to their larger size and their different trip characteristics. Therefore, a step-by-step approach is required to increase the feasibility and make proper investments for the initial deployment of electric trucks in the current systems.

Recent research papers have suggested optimization problems for truck routing with battery recharging/swapping infrastructure to overcome the low driving range of battery electric trucks. *(6,7)* have supported the deployment of electric trucks for delivery, which is usually known as electric vehicle routing problem (EVRP) or location and vehicle routing problem (LVRP). Mathematical models of the EVRPs find the best routing and recharging decision under the assumption that every electric truck can be recharged at the public truck recharging/swapping stations during its operation. Variables for purchasing battery electric truck, locating recharging/swapping stations, and routing each truck with enroute recharging are optimized to minimize the total cost, including purchasing, energy, and daily operational costs *(8)*.

The number of studies on the EVRP increased recently *(8)*, however, truck routing problems with enroute recharging/swapping assumptions still have limitations for real-world last-mile deliveries. First, the number of truck recharging stations in the real-world is too few to support enroute recharging. It is also difficult to increase the number of recharging stations in the near future because of the expensive purchasing cost for fast chargers (level 3), grid constraints, and difficulties to locate truck recharging stations in urban areas. Second, battery swapping technology is still in the early development stage, therefore, it needs various verification steps and improvements before it can be applied to real-world systems. Finally, fleet owners are hesitant to implement enroute recharging/swapping because it may increase the total operating time of drivers, which is strictly regulated.

Although freight truck electrification requires government subsidies, legal support, and technological improvements, the high energy efficiency of battery electric vehicles (BEVs) could reduce the burden of carriers in the long-term. Fleet electrification can also be an effective way to increase the number of locations allowing off-hour delivery (OHD) in urban areas since electric trucks offer much quieter and cleaner operations. OHD has been proposed as a way to reduce vehicle hours traveled (VHT) and traffic congestion by shifting some daytime delivery operations to nighttime. However, nighttime deliveries have a limitation due to emission, noise, and constraints on curb space *(9–11)*. Therefore, we suggest two-shift operations (daytime and nighttime deliveries) of electric trucks for the simultaneous implementation of both truck electrification and OHD in the e-commerce last-mile delivery. Nighttime operation is limited only to visiting locations that accept OHD, and every electric truck which has two routes per day should be recharged during break time between shifts at the depot. This multi-shift operation of electric trucks can be a low-hanging fruit for both fleet electrification and OHD implementations by overcoming the biggest hurdles for deployment: lack of public recharging facilities and low acceptance rate for OHD. This study will also make up for the current inefficiency of electric truck operations with lower driving range and higher purchasing cost, which will encourage fleet owners to have electric trucks for their last-mile deliveries.

In this paper, we propose an optimization model for the fleet-mix electric vehicle routing problem with multi-shift operations. The model (a) finds the best combination of trucks with heterogeneous engine types, (b) decides the best routes to deliver packages to several demand locations, (c) schedules two-shift operations for electric trucks including recharging between daytime and nighttime routes, and (d) allocates the proper number of EVSEs at the depot. The mathematical problem is described as mixed-integer linear programming (MILP) model to minimize the total operational cost, including total cost of ownership (TCO) of trucks, purchasing cost for EVSEs at the depot, and daily transportation cost. Also, a bi-level heuristic is suggested for efficient solution search, which consists of variable neighborhood search (VNS) at the upper-level and tabu search (TS) at the lower-level. The VNS heuristic finds the best combination of electric and conventional trucks by assigning pickup and delivery locations to each truck, while TS heuristic computes the shortest route to visit all the allocated locations of every truck and feeds the objective function value to the upper-level part.

All the solution algorithms are implemented in the POLARIS framework and evaluated with realistic simulations in Austin, Texas metropolitan area. POLARIS is an agent- and activity-based transportation simulation framework that provides, synthesizes, and simulates daily travel demands on the regional road network with links, traffic signs, and household locations *(12)*. We can get travel time and energy cost for every origin–destination pair at both daytime and nighttime and simulate the truck tours with background traffic on the POLARIS network. Based on synthesized e-commerce demand *(13),* the proposed algorithm computes the optimal fleet composition, including the EV ratio, the number of EVSEs, and daily routes, with given OHD acceptance ratio. We also analyze the cost efficiency, implementation feasibility, and the additional impact on the regional traffic network including total VHT, vehicle hours traveled (VMT), and traffic congestion.

The rest of this paper is organized as follows: related literature is reviewed in Section 2. The notation, assumptions, and MILP model are described in Section 3. In Section 4, the two-stage VNS-TS heuristic is proposed. Numerical experiments and simulated results are reported in Section 5, and conclusion and future studies are suggested in Section 6.

## 2. LITERATURE REVIEW

### 2.1. Electric Vehicle Routing Problems

Studies of the EVRP, LVRP for electric trucks, and green vehicle routing problem (GVRP) have been extended to cover cutting-edge technologies such as rapid charging *(7, 14, and 15),* battery swapping *(16–20),* and the powertrain mode selection for hybrid electric vehicles *(21).* Realistic representation of parameters and variables for non-linear recharging time *(22–24)*, waiting time at stations *(25–28]),* partial recharging strategy *(7, 29–31),* and energy consumption function during operations *(24, 32–35)* have also been suggested. Other research papers in *(31, 32, 36–39)* have considered mixed-fleet operation, including diesel-based and electrified trucks, and suggested mathematical problems to optimize the electrification ratio with minimum purchasing cost and operation cost. Common variables of proposed mathematical problems in this paragraph are: (a) locating recharging facilities, (b) finding routes with minimum travel time, and (c) recharging/swapping at public stations during operation.

Other research topics related to fleet electrification can be found in the strategic planning of EVSE allocation. Studies in *(40–46)* aim to locate recharging stations to maximize demand coverage, optimize the utilization of deployed recharging facilities, and suggest recharging strategies. Since these studies have usually not considered routing decisions and only focused on light-duty trucks or long-haul transport, it is limited for application to urban e-commerce delivery systems.

### 2.2. Routing Problems for Nighttime Deliveries

Only a few studies *(47–49)* have developed routing decisions for OHD implementation to analyze the traffic impact. Truck tours were synthesized using simple and constructive heuristic methods, and some routes were allocated to deliver during nighttime using iterative and scoring-based algorithms. Other related research topics can be found in multi-shift vehicle routing problems (VRPs) and pickup and delivery problems (PDPs) *(50–53)*. Because drayage operations are the target application of these research papers, long-haul transportation and full-truckload shipments are usually assumed. Therefore, it is also difficult to apply these models to the last-mile delivery with OHD.

The time window constraint is one of the common operational assumptions in VRPs and PDPs, which limits the arrival time at every delivery location to be within a feasible time range. This constraint can be used to represent a certain delivery location accepting OHD by setting the time window of the location to include nighttime. However, time window constraints usually make VRPs and PDPs much harder to solve since there are additional variables and inequalities that check the arrival time feasibility of all combinations of trucks and delivery locations *(54)*. Time window constraints were originally suggested to capture the passengers' inconvenience such as dial-a-ride services for patient transportation *(55)*, therefore, they may not be the best approach for modeling multi-shift operations, including nighttime deliveries. As a result, EVRPs with time windows *(6, 15, 28–31, 34–36)* do not consider the multi-shift operations and OHD. Finally, *(56)* studied the cost efficiency of OHD implementation in the urban business-to-business delivery system.

The most similar research has been suggested in *(57)*, which finds the best routes of electric trucks with depot recharging decisions. Every customer has multiple delivery locations, and each delivery location must be visited within a predetermined time window. The mathematical model finds the best routes of electric trucks to cover all customers' demand including recharging at depots. However, the model doesn't consider the break time and package sorting and loading during truck recharging at the depot, which represents every truck route as a single operation per day. To the best of our knowledge, no similar research topic has been reported in recent review papers *(8, 58, and 59)*.

## 2.3. Solution Search Methods

Recent studies on VRPs have focused on heuristic approaches because of their NP-hard nature *(54)*. Meta-heuristic algorithms are preferred to find near-optimal solutions for PDPs in a relatively short computational time, as reported by *(58, 59)*. We refer the reader to recent research papers *(15-17, 19–39)* for an extensive review about meta-heuristic approaches in PDPs and EVRPs: adaptive large neighborhood search (ALNS), tabu search (TS), genetic algorithm (GA), variable neighborhood search (VNS), iterative local search (ILS), and other hybrid meta-heuristics.

It is more challenging to implement meta-heuristics when heterogeneous fleets are assumed in VRPs and PDPs *(60)* because each truck type has different capabilities for routing and delivery. It is difficult to find better solutions with simple inter-route moves, and it is easy to be stuck in the local optima, therefore, more detailed procedures for assigning and/or swapping vehicles are required for mixed-fleet VRPs and PDPs. Inter-route reconstruction in ILS *(31)* and destroying and repairing procedures in ALNS or VNS *(31, 32, 36, 37, 60)* have been suggested to capture the heterogeneity of truck types.

## 3. PROBLEM DESCRIPTION

In this section, we introduce assumptions and the notation used throughout the paper. Consequently, we present a mathematical model of the electric vehicle routing problem with off-hour delivery (EVRP-OHD).

### 3.1. Assumptions

The EVRP-OHD assumes that every fleet owner can operate both conventional and electric trucks. Every diesel truck can have significantly longer route than an electrified truck to capture the current limitation on the driving range of electric trucks. On the other hand, each electrified truck can have two routes per day, visit demand locations allowing OHD at night, and recharge at its depot during the break time to overcome shorter driving range.

**Figure 1 (a and b)** are illustration of truck tours on a two-dimensional network with different truck types; a conventional truck can visit more demand locations with a single tour in **Figure 1 (a)** and an electric truck can visit all locations along two routes at both daytime and nighttime in **Figure 1 (b)**.

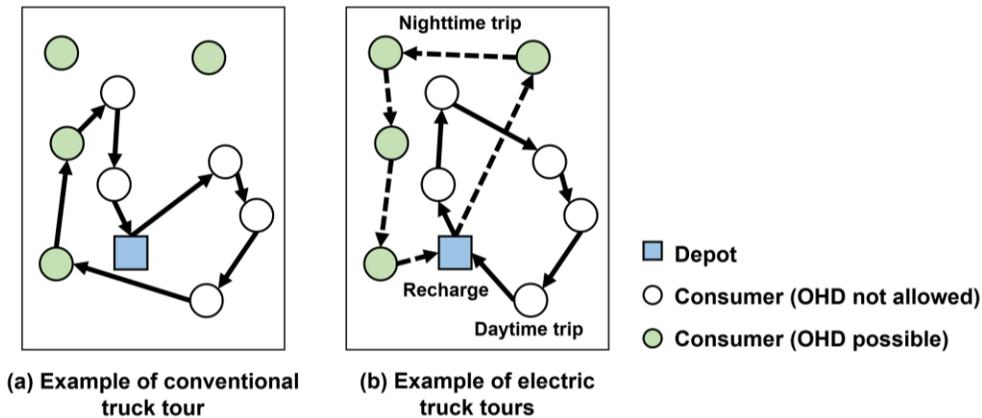

**Figure 1. Illustrations of the EVRP-OHD**

Although the freight truck electrification may enable fleet owners to reduce the number of trucks, sufficient purchasing of EVSEs at the depot is required for the recharging between daytime and nighttime operations. A static recharging assumption is proposed to compute the minimum number of EVSEs. As a result, fleet electrification has a trade-off between reducing the number of trucks and the high purchasing costs of both electric trucks and EVSEs.

We apply common assumptions and constraints for VRPs, such as (a) every demand location must be visited by a single truck, (b) each truck tour should operate within maximum driving hours and loading capacity, and (c) each tour must be a closed loop.

## 3.2. Mathematical Formulation

We define the EVRP-OHD on a directed graph $\mathcal{G} = (\mathcal{N}, \mathcal{A})$ with location set $\mathcal{N} = \mathcal{N}^C \cup \mathcal{N}^D \cup \mathcal{N}^O$ and arc set $\mathcal{A}$. Subset $\mathcal{N}^C$ represents the depot location of every fleet owner. Two subsets for customer locations are suggested, $\mathcal{N}^D$ is a subset of demand locations not allowing nighttime delivery, and $\mathcal{N}^O$ consists of locations allowing OHD. We also assume that all locations in $\mathcal{N}^O$ can be visited by electrified trucks only because they produce far less noise. Let $\mathcal{K}$ denote the set of trucks, and two types of trucks are defined: conventional trucks ($k \in \mathcal{K}^C$) and electrified trucks ($k \in \mathcal{K}^E$).

We assume different purchasing costs for conventional ($C^C$), electric trucks ($C^E$), and EVSEs ($C^{EVSE}$). Travel time from location $i$ to $j$ during daytime ($T_{ij}^D$) and nighttime ($T_{ij}^O$) are obtained through POLARIS simulation. E-commerce demand at location $i$ ($V_i$) is estimated based on the model developed in *(13)*, where the e-commerce model attributes are obtained though the population synthesizer of POLARIS. Drop-off time at location $i$ ($S_i$) is proportional to the demand volume at the location, as $S_i = V_i \times$ unit drop-off time. Heterogeneous truck capacity ($Q^C$ and $Q^E$) and driving range ($D^C$ and $D^E$) are considered for conventional trucks and electric trucks respectively. $H^{EVSE}$ is unit recharging speed, $B$ is break time between daytime and nighttime operations, and $bigM$ is a big number.

Let $u_i$ be the number of EVSEs at the depot of carrier establishment $i$ for every $i \in \mathcal{N}^C$, and variable $r_{ik}$ is recharging time for every $k \in \mathcal{K}^E$ and $i \in \mathcal{N}^C$. The binary variable $w_k$ indicates whether the truck $k$ is assigned for the operation ($w_k = 1$) or not ($w_k = 0$) for every truck $k \in \mathcal{K}$. Finally, variables $x_{ijk}$ and $y_{ijk}$ be a binary variable equal to 1 if and only if the truck $k$ visits arc $(i,j)$ during daytime and nighttime, respectively, for every $(i,j) \in \mathcal{A}$ and $k \in \mathcal{K}$.

**Table 1** shows sets, parameters, and variables in the EVRP-OHD.

**Table 1. Notations**

| **Sets** |
| --- |
| $\mathcal{N}$: the set of locations ($\mathcal{N} = \mathcal{N}^C \cup \mathcal{N}^D \cup \mathcal{N}^O$) |
| $\mathcal{N}^C$: the subset of depots |
| $\mathcal{N}^D$: the subset of customer locations not allowing off-hour delivery |
| $\mathcal{N}^O$: the subset of customer locations allowing off-hour delivery |
| $\mathcal{A}$: the set of arcs |
| $\mathcal{K}$: the set of trucks ($\mathcal{K} = \mathcal{K}^C \cup \mathcal{K}^E$) |
| $\mathcal{K}^C$: the subset of conventional trucks |
| $\mathcal{K}^E$: the subset of electrified trucks |
| **Parameters** |
| $C^C$: the TCO of a conventional truck |
| $C^E$: the TCO of an electrified truck |
| $C^{EVSE}$: the purchasing cost of an EVSE |
| $T_{ij}^D$: travel time from location $i$ to $j$ during daytime |
| $T_{ij}^O$: travel time from location $i$ to $j$ during nighttime |
| $Q^C$: transportation capacity of a conventional truck |
| $Q^E$: transportation capacity of an electrified truck |
| $D^C$: driving range of a conventional truck |
| $D^E$: driving range of an electrified truck |
| $V_i$: demand volume at the location $i$ |
| $S_i$: service time for drop-off at the location $i$ |
| $H^{EVSE}$: unit recharging speed (traveling time / recharging time) |
| $B$: break time between daytime and off-hour deliveries |
| $bigM$: a big number |
| **Variables** |
| $w_k$: = 1 if the truck $k$ is assigned for delivery, = 0 otherwise |
| $x_{ijk}$: = 1 if the truck $k$ departs from the location $i$ to $j$ at daytime, = 0 otherwise |
| $y_{ijk}$: = 1 if the truck $k$ departs from the location $i$ to $j$ at nighttime, = 0 otherwise |
| $z_{ik}$: total volume of packages when the truck $k$ at the location $i$ |
| $u_i$: the number of EVSEs at the location $i$ |
| $r_{ik}$: the total recharging time of the electrified truck $k$ during break time at the depot $i$ |

Every carrier establishment minimizes the total cost, which is the sum of TCO of both conventional and electrified trucks, EVSE purchasing cost, and daily operation cost. The EVRP-OHD can be formulated as the multi-objective MILP model, as follows:

$$\text{Minimize} \quad \sum_{k \in \mathcal{K}^C} C^C w_k + \sum_{k \in \mathcal{K}^E} C^E w_k + \sum_{i \in \mathcal{N}^C} C^{EVSE} u_i \quad (1\text{-}1)$$

$$\text{Minimize} \quad \sum_{k \in \mathcal{K}^E} \sum_{i,j \in \mathcal{N}} \left( T_{ij}^D x_{ijk} + T_{ij}^O y_{ijk} \right) \quad (1\text{-}2)$$

subject to

**Day-time deliveries**

$$\sum_{j \in \mathcal{N}} x_{ijk} \leq w_k \quad \forall i \in \mathcal{N}^C, k \in \mathcal{K} \quad (2)$$

$$\sum_{j \in \mathcal{N}} x_{ijk} = \sum_{j \in \mathcal{N}} x_{jik} \quad \forall i \in \mathcal{N}, k \in \mathcal{K} \quad (3)$$

$$\sum_{j \in \mathcal{N}} \sum_{k \in \mathcal{K}} x_{ijk} = 1 \quad \forall i \in \mathcal{N}^D \quad (4)$$

$$\sum_{i,j \in \mathcal{N}} T_{ij}^D x_{ijk} \leq D^C \quad \forall k \in \mathcal{K}^C \quad (5)$$

$$z_{jk} \geq z_{ik} + V_j - bigM(1 - x_{ijk}) \quad \forall i \in \mathcal{N}, j \in \mathcal{N}^C, k \in \mathcal{K} \quad (6)$$

$$z_{ik} \leq Q^C \quad \forall i \in \mathcal{N}, k \in \mathcal{K}^C \quad (7)$$

**Off-hour deliveries**

$$\sum_{j \in \mathcal{N}} y_{ijk} \leq \sum_{j \in \mathcal{N}} x_{ijk} \quad \forall i \in \mathcal{N}^C, k \in \mathcal{K}^E \quad (8)$$

$$\sum_{j \in \mathcal{N}} y_{ijk} = \sum_{j \in \mathcal{N}} y_{jik} \quad \forall i \in \mathcal{N}^O, k \in \mathcal{K}^E \quad (9)$$

$$\sum_{j \in \mathcal{N}} \left( \sum_{k \in \mathcal{K}} x_{ijk} + \sum_{k \in \mathcal{K}^E} y_{ijk} \right) = 1 \quad \forall i \in \mathcal{N}^O \quad (10)$$

$$\sum_{i,j \in \mathcal{N}} T_{ij}^D x_{ijk} \leq D^E \quad \forall k \in \mathcal{K}^E \quad (11)$$

$$\sum_{i,j \in \mathcal{N}} T_{ij}^O y_{ijk} \leq D^E \quad \forall k \in \mathcal{K}^E \quad (12)$$

$$z_{jk} \geq z_{ik} + V_j - bigM(1 - y_{ijk}) \quad \forall i \in \mathcal{N}, j \in \mathcal{N}^O, k \in \mathcal{K}^E \quad (13)$$

$$z_{ik} \leq Q^E \quad \forall i \in \mathcal{N}, k \in \mathcal{K}^E \quad (14)$$

**Recharging**

$$H^{EVSE} r_{ik} \geq \sum_{i,j \in \mathcal{N}} C_{ij}^E x_{ijk} - bigM \left( 1 - \sum_{j \in \mathcal{N}^O} y_{ijk} \right) \quad \forall i \in \mathcal{N}^C, k \in \mathcal{K}^E \quad (15)$$

$$\sum_{k \in \mathcal{K}} r_{ik} \leq B u_i \quad \forall i \in \mathcal{N}^C, k \in \mathcal{K}^E \quad (16)$$

**Variable constraints**

$$u_i \in \mathbb{N}^+ \quad \forall i \in \mathcal{N}^C$$
$$w_k, x_{ijk}, y_{ijk} \in \{0,1\} \quad \forall i,j \in \mathcal{N}, k \in \mathcal{K} \quad (17)$$
$$r_{ik}, z_{ik} \in \mathbb{R}^+ \quad \forall i \in \mathcal{N}, k \in \mathcal{K}$$

Objective function (1-1) minimizes the sum of TCO of trucks and EVSE purchasing cost and (1-2) optimizes each truck's daily route with the minimum total travel time. Constraints (2) and (3) represent the route construction constraints for every assigned truck. Every demand location must be visited by one of the trucks by constraint (4). Constraints (5) and (6–7) are truck capacity and maximum driving range constraints, respectively. Constraints (8–14) are only for OHDs, which shares similar structure with constraints (2–7). Constraint (15) ensures that every electrified truck requires recharging during break time if and only if the truck has deliveries at nighttime. Constraint (16) makes the sum of recharging time be less than break time by assigning proper number of EVSEs. Finally, constraint (17) is a variable constraint.

## 4. BI-LEVEL VNS-TS HEURISTIC

We developed a bi-level VNS-TS heuristic; truck assignment in the upper-level part and route improvement in the lower-level part. Truck assignment, engine type, and visiting locations vector are updated with VNS heuristic in the upper-level, and the lower-level problem finds the best routes of updated trucks with TS procedures. **Figure 2** below shows the structure of the VNS-TS heuristic algorithm.

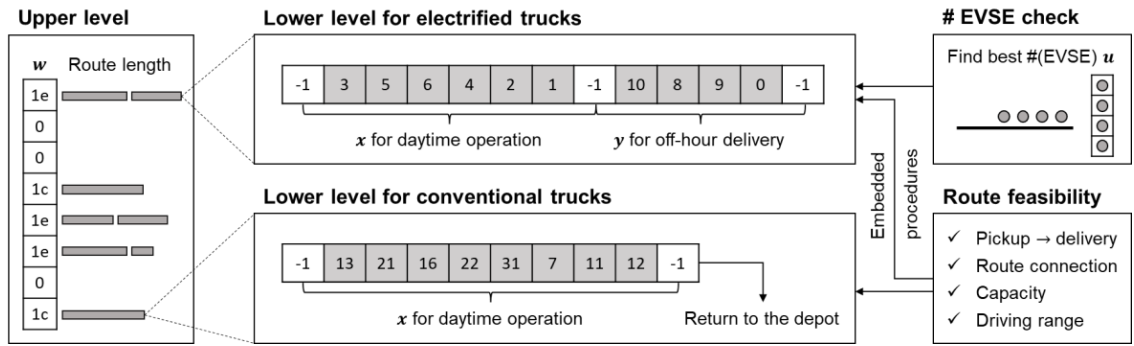

**Figure 2. Bi-level VNS-TS heuristic description**

In the upper-level part of **Figure 2**, there is a vector which shows the assignment and type info of each truck; '0' for non-assigned, '1c' for assigned conventional, and '1e' for assigned electrified trucks. Also, it has a location vector to be visited by each assigned truck for the demand covering constraint. Upper-level solutions are updated using VNS heuristic with several inter-route movements to find the best fleet mix to cover all the demands with the minimum purchasing costs (1-1). Whenever upper-level problem needs to compute objective function of a certain solution, lower-level problem is called to find the best order of visiting locations of a truck and return the value of the objective function (1-2) to the upper-level. TS heuristic updates string-structured solutions with intra-route movements. There are two embedded procedures to compute the EVSE purchasing cost and guarantee the feasibility of lower-level solutions.

The bi-level structure offers efficient solution search for large-sized instances of the EVRP-OHD, encompassing the entire e-commerce delivery system in simulated metropolitan cities, by partitioning the solution space through variable classification. While the bi-level heuristic efficiently manages the number of neighborhood solution, it generally requires more iterations to converge. Further, the size of EVRP-OHD instances is significantly larger than that of typical VRPs, this paper focuses on heuristic development that emphasizes exploration at the upper level with ample randomness, while simplifying the lower level to minimize overall computational demands.

All the parts for both upper-level and lower-level problem are described below.

### 4.1. VNS for Upper-level Problem

The VNS heuristic explores a local optimum with a local search method with several neighborhood structures and perturbs the current solution to escape from the corresponding valley in every iteration. The main local searching steps of the VNS are (1) search more distant neighborhoods of the current solution when it fails to find any improved solutions, (2) jump to the improved solution, and (3) repeat the local search with adjacent neighborhoods. The heuristic also has a shaking procedure to add randomness. Algorithm 1 is the pseudo-code of the upper-level VNS heuristic.

---

**Algorithm 1.** Variable Neighborhood Search

**start**
    $sol = InitialSolutionByClustering(|\mathcal{K}|)$
    **repeat**
        **for** $n$ **in** $NeighborhoodStructure$
            $sol' = Shake(sol, n)$
            $SOL_n = NeighborhoodSet(SOL', n)$
            $sol^* = \underset{s \in SOL_n}{\mathrm{argmin}}\{ObjValue(s) \mid Feasible(s) = true\}$
            **if** $ObjValue(sol^*) < ObjValue(sol)$
                $sol = sol^*$
                $n = 1$
            **else**
                $n = n + 1$
            **end if**
        **end for**
    **return** $sol$ **if** $time > MaxComputationTime$
**end**

---

We implement the VNS with a simple initialization procedure, which generates clusters of customers by using k-medoid algorithm. The number of clusters is set to be the minimum number of diesel trucks to cover the whole customer volume ($k = |\mathcal{K}^C|$). k medoids are selected randomly among customers, all other customers are assigned to the closest medoid, and all constraints in the EVRP-OHD are not considered in this step. Finally, a single diesel truck is assigned to every cluster, and the lower-level heuristic completes the initial tour of each truck.

Common shaking procedures randomly choose a neighborhood solution to be a current solution at every iteration in previous works *(60)*, but it is not enough to add sufficient randomness to the problem when the number of locations is greater than several thousands. Therefore, shaking procedure in this paper prioritizes trucks randomly, and the neighborhood solution preferentially changes solution structure with highly prioritized truck. It enables not only to reduce the number of neighborhood solutions efficiently but also to add randomness for preventing solutions stuck in the local optima.

The initial solution is updated with several neighborhood structures, and **Figure 3** below is an illustration of the neighborhood generation. Since we sorted neighborhood structures from the smallest to the biggest changes in the current solution as '1–insert', '1–swap', '2–insert', 'type change', 'route insert' and 'route delete', the heuristic explores the solution space farther when the previous neighborhood structure cannot find better solutions. Whenever the location subset of a certain truck is changed, the lower-level heuristic is called to find the best route and to return the total travel time of the truck. All lower-level procedures are parallelized with agent-based modeling framework, which enables efficient solution evaluation.

Current solution is always updated whenever better feasible solutions have found, and the whole procedure is repeated until the termination criteria met.

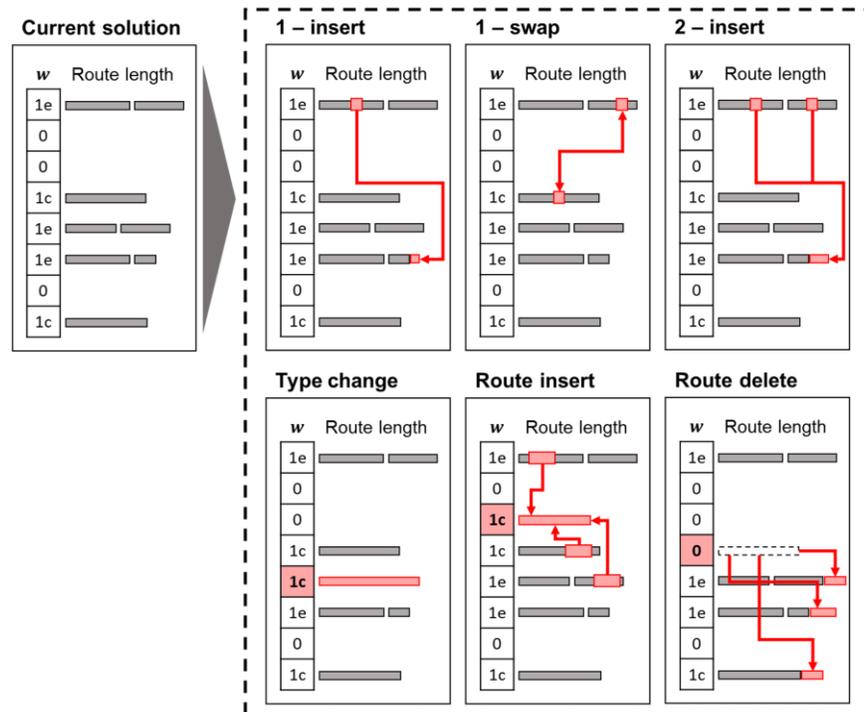

**Figure 3. An illustration for neighborhood structure in VNS heuristic**

## 4.2. TS for Lower-level Problem

The TS heuristic uses a trajectory-based search method. The basic idea of the TS heuristic is moving a current solution to one of the locally improved solutions and recording some elements of a solution movement in the tabu list. The elements in the list cannot be revisited for the next $\lambda$ iterations, denoted as the tabu tenure. This process prevents repetition in a sequence of solutions, which may reduce the frequency of getting stuck in the local optima. The TS heuristic finds the best solution until its termination, which is usually the maximum iterations at local optima or the maximum computation time.

In this paper, common features and procedures of the TS heuristic are applied for the lower-level problem. String-structured solution is suggested in **Figure 2**, and we add one additional node for the recharging during breaktime when the truck is electrified. Neighborhood solutions are generated with 'node-swap', 'insert', and '2-opt' procedures as shown in **Figure 4**. To prevent repeated searches, we applied the common structure of tabu-list from *(25,30)* – containing moved location ids – with size = $\sqrt{|\mathcal{L}|}$ where $|\mathcal{L}|$ is the length of the solution string. Because every lower-level problem only solves a single-truck EVRP with a few locations already assigned to the truck, only few iterations ($=\sqrt{2 \times |\mathcal{L}|}$) are allowed to reduce the computational requirements at the lower-level.

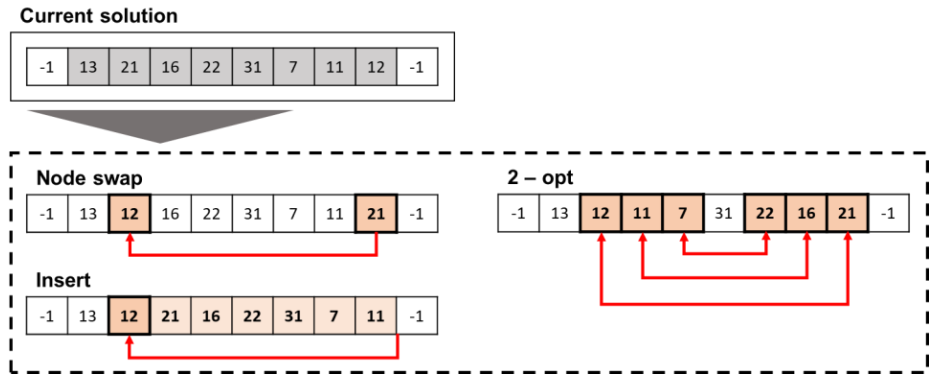

**Figure 4. An illustration for neighborhood structure in TS heuristic**

## 5. EXPERIMENTAL RESULTS

In this section, the performance of the suggested bi-level VNS-TS heuristic is evaluated with simulated e-commerce demand of Bloomington, Illinois in POLARIS and compared to the commercial optimization solver. Also, the EVRP-OHD is embedded in the POLARIS framework and simulated to represent the entire residential e-commerce deliveries in the metropolitan area of Austin, Texas within POLARIS. The efficiency of the joint implementation of fleet electrification and off-hour delivery operation is analyzed and reported for several scenarios.

## 5.1. Performance Evaluation

A few e-commerce orders in Bloomington were sampled including location-to-location travel time, number of packages, and drop-off time. Also, locations allowing nighttime delivery are randomly selected based on 4 OHD ratios. All other parameters related to truck TCO, EVSE purchasing cost, truck capacity and operation time limitation were simplified to be homogeneous, regardless of truck type. Finally, the proposed bi-level VNS-TS heuristic was compared with IBM ILOG CPLEX 12.10 on the same personal computer (Intel Xeon W-2155 @ 3.30GHz, 20 threads, and 32GB ram).

Simplified parameters are assumed for clearer comparison of performance between CPLEX solver and the suggested heuristic. The TCO, driving range, and delivery capacity of electric trucks are assumed to be same to those of diesel trucks. Also, EVSE purchasing cost and recharging time are not considered to simplify the recharging decisions. Customer locations accepting OHD are selected randomly with given acceptance ratio and average travel speed during nighttime is 30% higher than that during daytime to force trucks operate during nighttime. Finally, the MILP model and the VNS-TS heuristic minimize the sum of objective functions (1-1) and (1-2) with artificial TCO value ($C^C = C^E = 1,000$).

**Table 2** shows results of the sampled Bloomington e-commerce deliveries from both CPLEX and the VNS-TS heuristic. The number of sample locations ($|\mathcal{N}|$), the fleet size ($|\mathcal{K}|$), and the OHD acceptance ratio are reported in the scenario column. The objective function value and computation time are reported in each solution approach section. Since CPLEX is allowed to search solution within 1 hour only, the gap ratio between the best feasible solution and the lower bound is reported as MIP gap. Finally, the difference of the objective function value between CPLEX and the VNS-TS heuristic is reported in the last column.

Bi-level heuristic could report better solution than that from CPLEX in almost every scenario. The total travel time can be reduced with more deliveries during nighttime, and these trends – higher solution value with fewer customers allowing OHD – can be easily observed with the results from the VNS-TS heuristic in **Table 2.** However, CPLEX usually returned higher solution value even though more customers accept delivery during nighttime. MIP gap also tends to increase with higher OHD ratio, and these computational inefficiencies are because the mathematical model needs to search larger space of decision variable $y_{ijk}$ ($|Y| \propto |\mathcal{N}|^2 \times \text{OHDratio}^2$).

In general, all scenario tests cannot be finished within a given time limit of 1 hour. CPLEX reported poor gaps between feasible solutions and lower bounds, which show the complexity of the EVRP-OHD. We believe that the EVRP-OHD is much harder to find the optimal solution than the common VRPs because the model has (1) two sets of three-index decision variables for routing at both daytime ($x_{ijk}$) and nighttime ($y_{ijk}$), (2) objective function with weighted sum of both strategic purchasing cost and daily operational cost, and (3) more constraints with big-M notations to indicate the condition about multi-shift operation and recharging.

**Table 2. Result comparison between CPLEX and VNS-TS heuristic**

| Scenario | | Commercial solver (CPLEX) | | | VNS-TS heuristic | | gap B-A (%) |
|---|---|---|---|---|---|---|---|
| $|\mathcal{N}|, |\mathcal{K}|$ | OHD Ratio (%) | Obj. value (A) | MIP gap (%) | Computation time (sec) | Obj. value (B) | Computation time (sec) | |
| {20, 4} | 0.0 | 20,761 | 85.15 | 3,600 | 20,812 | 1 | 0.25 |
| | 20.0 | 20,762 | 82.74 | 3,600 | 20,265 | < 1 | -2.40 |
| | 40.0 | 20,762 | 82.77 | 3,600 | 18,424 | < 1 | -11.26 |
| | 100.0 | 20,959 | 84.36 | 3,600 | 14,539 | < 1 | -30.63 |
| {33, 10} | 0.0 | 39,183 | 86.50 | 3,600 | 38,402 | 1 | -1.99 |
| | 20.0 | 40,521 | 84.42 | 3,600 | 35,442 | 1 | -12.54 |
| | 40.0 | NONE | - | 3,600 | 33,139 | 1 | - |
| | 100.0 | NONE | - | 3,600 | 27,973 | 1 | - |
| {37, 8} | 0.0 | 33,306 | 79.81 | 3,600 | 31,423 | 1 | -5.65 |
| | 20.0 | 34,490 | 87.00 | 3,600 | 31,071 | 13 | -9.91 |
| | 40.0 | 34,008 | 86.91 | 3,600 | 30,173 | 2 | -11.28 |
| | 100.0 | NONE | - | 3,600 | 25,257 | 2 | - |
| {46, 10} | 0.0 | 40,663 | 86.11 | 3,600 | 39,676 | 5 | -2.43 |
| | 20.0 | 44,023 | 85.32 | 3,600 | 34,626 | 17 | -21.35 |
| | 40.0 | 39,906 | 85.18 | 3,600 | 36,527 | 65 | -8.47 |
| | 100.0 | NONE | - | 3,600 | 33,372 | 2 | - |
| {49, 8} | 0.0 | 38,730 | 87.45 | 3,600 | 35,703 | 3 | -7.81 |
| | 20.0 | 37,857 | 90.77 | 3,600 | 31,555 | 14 | -16.65 |
| | 40.0 | 36,662 | 85.34 | 3,600 | 32,056 | 4 | -12.56 |
| | 100.0 | NONE | - | 3,600 | 29,313 | 1 | - |
| {58, 8} | 0.0 | 38,865 | 87.39 | 3,600 | 35,590 | 8 | -8.43 |
| | 20.0 | 37,707 | 87.45 | 3,600 | 35,127 | 31 | -6.84 |
| | 40.0 | 37,619 | 86.58 | 3,600 | 33,662 | 71 | -10.52 |
| | 100.0 | NONE | - | 3,600 | 29,529 | 1 | - |

※ NONE: Feasible solutions cannot be found within the max CPU time

## 5.2. Scenario Evaluation

Demand locations and all other details for e-commerce delivery are synthesized using the e-commerce choice model based on the synthesized population within POLARIS. Below, **Figure 5** is the target study area in POLARIS, which represents the Austin metropolitan area in Texas with 39,638 aggregated household locations of residential land-uses, 15,833 links, and 2,161 zones. Also, the major four e-commerce companies are emulated in the POLARIS-Austin, operating out of 20 different depots. There are an average of 280,000 orders among 830,000 households. Each company partitions Austin into several sub-zone clusters and each depot serves customers in a cluster as shown in **Figure 5**. As a result, we implemented bi-level heuristic 20 times for solving single-depot EVRP-OHD.

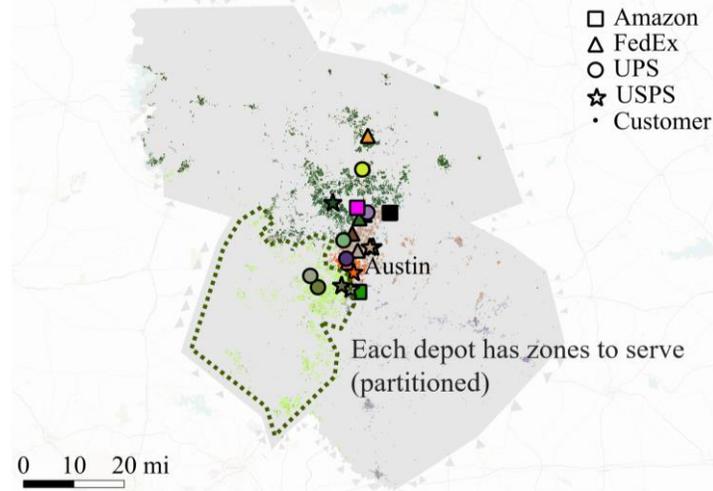
**Figure 5. Target network for numerical experiments**

The EVRP-OHD used travel time between every pair of origin and destination locations from POLARIS network with background traffic, which enables getting reasonable travel time during both daytime and nighttime. We considered 10 test scenarios with 2 EV ranges for delivery trucks and 5 cases of OHD acceptance ratio. The TCO of an electric truck and purchasing cost of a recharging plug is assumed to be 95% and 20% of the TCO of diesel truck respectively *(61)*. We set $R^E$ = 7 hours based on the current technology of MD electric truck by dividing 210 miles into 30 MPH and $R^E$ = 11 hours when EV range becomes 50% longer. Also, locations allowing nighttime delivery are randomly chosen with given OHD ratio. Other parameters are assumed for test experiments such as: $S_i/V_i$ (unit drop-off time) = 120 secs, $Q^C = Q^E$ = 140 orders, $R^C$ = 12 hours, $H^{EVSE}$ = 1 hour recharging after 7 hour driving, and $B$ = maximum 4 hours. Finally, the VNS-TS heuristic minimizes the total TCO and EVSE purchasing cost first at the upper-level and then finds the best tours at the lower-level because the actual value of strategic purchasing cost in objective function (1-1) and the daily travel time in (1-2) are not comparable.

Due to the computational intensity of repeatedly solving the single-depot EVRP-OHD problem with many thousands of locations, we conducted all simulations on the Broadwell partition of the BEBOP cluster operated by Laboratory Computing Resource Center at Argonne National Laboratory. The cluster consists of over 650 high-performance compute nodes, each with 36 physical cores and 128GB DDR4 RAM. By integrating our proposed methodology within the highly efficient POLARIS framework, we were able to run all simulations for solving EVRP-OHD using between 5 (OHD=0%) and 12 (OHD=50%) node hours per simulation.

**Table 3** shows results of ten scenarios in POLARIS-Austin.

**Table 3. Test results**

| Scenario | | Fleet composition | | | Simulation results | | | | | |
|---|---|---|---|---|---|---|---|---|---|---|
| EV range (hours) | OHD ratio (%) | Trucks (A) | Routes (B) | EV ratio (B-A)/A, (%) | VMTs (miles) | VHTs (hours) | Avg MPH | Avg VMT (miles/B) | Avg VHT (hours/B) | Avg OH (hours/A) |
| 7.0 | 0.0 | 2,210 | 2,210 | 0.0 | 60,847.7 | 2,032.4 | 29.9 | 110.1 | 3.7 | 7.9 |
| 7.0 | 10.0 | 2,531 | 2,980 | 17.7 | 68,989.8 | 2,306.3 | 29.9 | 92.6 | 3.1 | 7.4 |
| 7.0 | 20.0 | 2,471 | 3,373 | 36.5 | 71,263.0 | 2,373.6 | 30.0 | 84.6 | 2.8 | 7.7 |
| 7.0 | 30.0 | 2,354 | 3,650 | 55.1 | 73,700.3 | 2,428.2 | 30.4 | 80.8 | 2.7 | 8.1 |
| 7.0 | 50.0 | 2,186 | 3,992 | 82.7 | 76,604.9 | 2,493.1 | 30.7 | 76.8 | 2.5 | 8.9 |
| 11.0 | 0.0 | 2,221 | 2,221 | 0.0 | 61,066.6 | 2,051.6 | 29.8 | 110.0 | 3.7 | 7.9 |
| 11.0 | 10.0 | 2,240 | 2,828 | 26.3 | 65,799.5 | 2,192.9 | 30.0 | 93.1 | 3.1 | 8.1 |
| 11.0 | 20.0 | 2,045 | 3,051 | 49.2 | 68,679.3 | 2,254.1 | 30.5 | 90.1 | 3.0 | 9.0 |
| 11.0 | 30.0 | 1,856 | 3,107 | 67.4 | 69,934.9 | 2,286.9 | 30.6 | 90.1 | 2.9 | 10.0 |
| 11.0 | 50.0 | 1,734 | 3,086 | 78.1 | 69,973.4 | 2,250.8 | 31.1 | 90.8 | 2.9 | 10.6 |

※ Average number of demand packages = 282,257    ※ OH: operation hours = total travel time + total drop-off time

As shown in **Figure 6 (a)** and **(b)**, a higher household ratio accepting OHD requires more electric trucks and increases delivery stops during nighttime. These values monotonically increase with higher OHD acceptance ratio whether the EV driving range is short or long. Also, about 40% of the household locations accepting OHD are actually delivered during nighttime from the EVRP-OHD heuristic results.

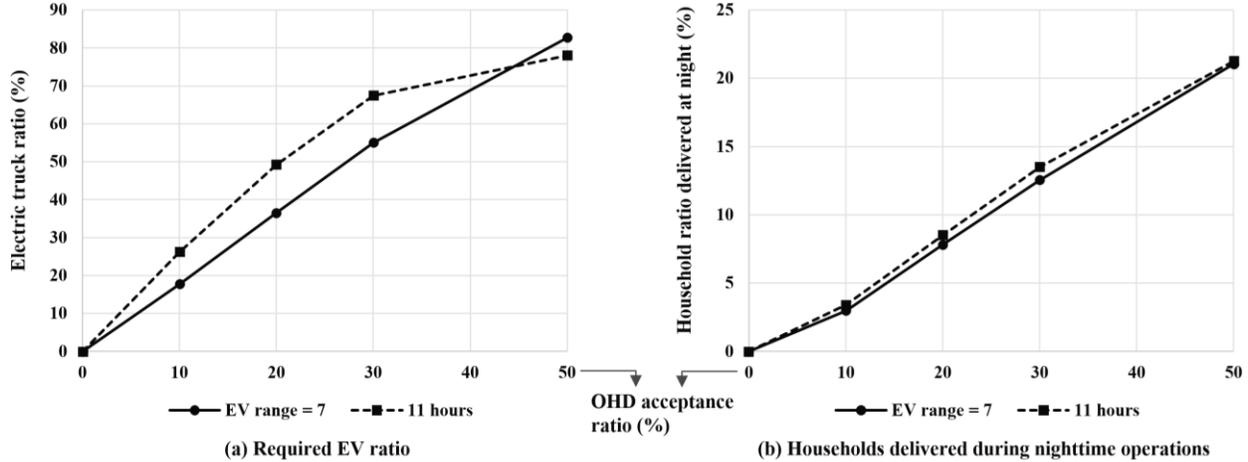

(a) Required EV ratio    (b) Households delivered during nighttime operations

**Figure 6. Fleet composition and size**

One of the benefits of OHD implementation is avoiding traffic by moving some delivery trips into nighttime. In **Table 3**, average truck speed increased with higher OHD acceptance ratio because of fewer traffic on the network, which also makes shorter average vehicle hour traveled (not including drop off times) as shown in **Figure 7 (a)**. Average VMTs in **Figure 7 (b)** have also been decreased with more frequent nighttime operations, which implies that electric trucks can not only avoid traffic congestion during daytime but also find shorter paths during nighttime. Since the average miles and hours decrease significantly with relatively low OHD acceptance ratios of 10–20%, joint implementation of fleet electrification and nighttime delivery to limited demand locations could be an initial approach to reduce energy consumption in real-world e-commerce delivery systems.

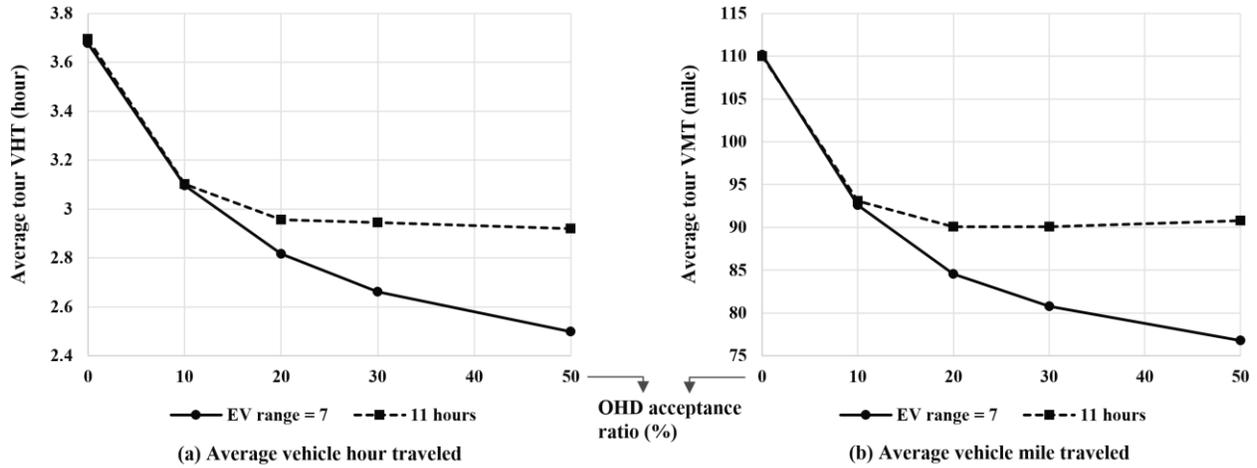

**Figure 7. Average truck speed and VHTs**

Although electric truck is preferred with higher OHD acceptance ratio, more trucks could be required. In **Figure 8 (a)**, we compared the overall fleet size, and shorter EV range scenario increased the fleet size when OHD acceptance ratio is 10–30%. On the other hand, long EV range could reduce the fleet size from 2,221 to 1,734. Also, **Figure 8 (b)** shows that average operation hours of every truck could increase because the electrified trucks need more trips departing from and returning to its depot for recharging. In **Table 3**, total number of tours, VMTs, and VHTs also increased with the same reason. Therefore, locating depots closer to demand locations is important to resolve these issues, and pairing the EVRP-OHD with micro-hub operations could be an option to reduce the overall traveling distance from depots to demand locations.

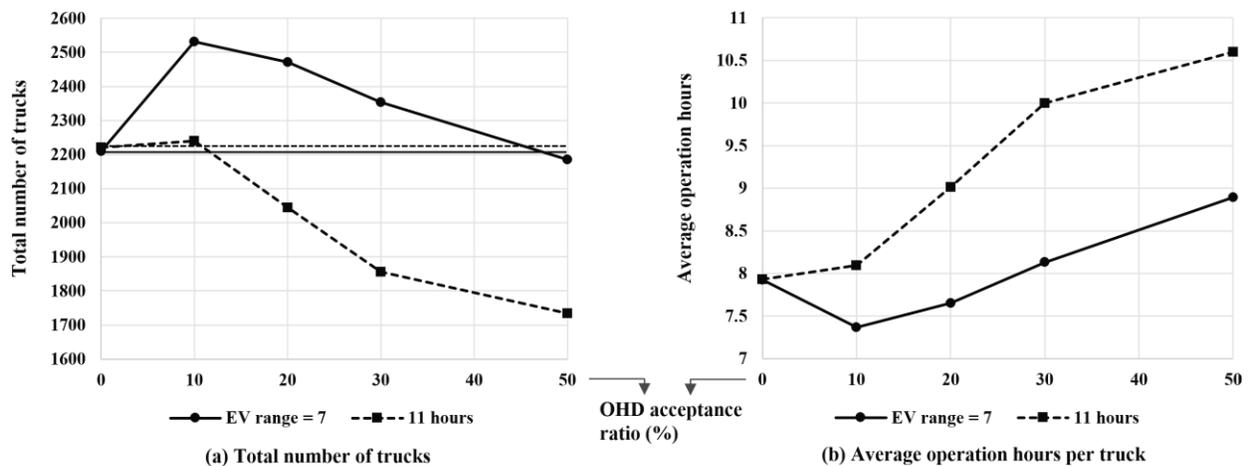

**Figure 8. Fleet composition and size**

# 6. CONCLUSION

In this paper, we proposed the fleet-mix electric vehicle routing problem with limited off-hour delivery implementation (EVRP-OHD) to initiate fleet electrification decisions in e-commerce delivery systems. The EVRP-OHD considered the fleet electrification under current limitation of EV technologies, such as lack of public recharging stations, shorter driving ranges than conventional diesel trucks, and more expensive purchasing costs for EV and EVSE. We also applied limited off-hour delivery implementation to the model based on the strengths of electric trucks, including clean and quiet operations. Every electric truck can have more flexible routes and visit demand locations accepting OHD during nighttime delivery. We suggested a mathematical formulation to minimize both purchasing cost of trucks and EVSEs and daily operation costs and designed a bi-level VNS-TS heuristic for efficient solution search.

Test experiments with several EV driving ranges and OHD acceptance ratio scenarios were performed in the POLARIS framework. The results demonstrated the efficiency of deploying freight electric trucks for OHD implementations. In general, a low OHD acceptance ratio (10–20% of total demand) was enough to get the benefit of multi-shift operations. Applying the EVRP-OHD led to decreases of more than 20% of average VMT and VHT, indicating that this may also reduce energy usage for e-commerce delivery operations. Also, the EV ratio required to support the optimal strategy also increased with higher OHD acceptance ratio even though we applied limited operational ranges for electric trucks. Results were highly dependent on the driving range of EVs, therefore, fleet electrification can be boosted with higher EV ranges.

In this study, we limited nighttime deliveries to being operated by electric trucks only. In the future this constraint could be relaxed allow conventional truck types to deliver overnight where municipality policies accept nighttime operations. Also, the EVRP-OHD could be extended to include delivery lockers, micro-hubs, and other technologies to improve drivers' safety during nighttime operations and reduce the overall traveling distance between depots to demand locations in future studies.

# ACKNOWLEDGMENTS

This report and the work described were sponsored by the U.S. Department of Energy (DOE) Vehicle Technologies Office (VTO) under the Systems and Modeling for Accelerated Research in Transportation (SMART) Mobility Laboratory Consortium, an initiative of the Energy Efficient Mobility Systems (EEMS) Program. The submitted manuscript has been created by the UChicago Argonne, LLC, Operator of Argonne National Laboratory (Argonne). Argonne, a U.S. Department of Energy Office of Science laboratory, is operated under Contract No. DE-AC02-06CH11357. The U.S. Government retains for itself, and others acting on its behalf, a paid-up nonexclusive, irrevocable worldwide license in said article to reproduce, prepare derivative works, distribute copies to the public, and perform publicly and display publicly, by or on behalf of the Government.

16. Yang, J., and Sun, H. Battery swap station location-routing problem with capacitated electric vehicles. *Computers and operations research*, 2015. 55:217–232.

17. Hof, J., Schneider, M., and Goeke, D. Solving the battery swap station location-routing problem with capacitated electric vehicles using an AVNS algorithm for vehicle-routing problems with intermediate stops. *Transportation research part B: methodological*, 2017. 97:102–112.

18. Wang, Y. W., Lin, C. C., and Lee, T. J. Electric vehicle tour planning. *Transportation Research Part D: Transport and Environment*, 2018. 63:121–136.

19. Masmoudi, M. A., Hosny, M., Demir, E., Genikomsakis, K. N., and Cheikhrouhou, N. The dial-a-ride problem with electric vehicles and battery swapping stations. *Transportation research part E: logistics and transportation review*, 2018. 118:392–420.

20. Jie, W., Yang, J., Zhang, M., and Huang, Y. The two-echelon capacitated electric vehicle routing problem with battery swapping stations: Formulation and efficient methodology. *European Journal of Operational Research*, 2019. 272(3):879–904.

21. Zhen, L., Xu, Z., Ma, C., and Xiao, L. Hybrid electric vehicle routing problem with mode selection. *International Journal of Production Research*, 2020. 58(2):562–576.

22. Montoya, A., Guéret, C., Mendoza, J. E., and Villegas, J. G. The electric vehicle routing problem with nonlinear charging function. *Transportation Research Part B: Methodological*, 2017. 103:87–110.

23. Froger, A., Mendoza, J. E., Jabali, O., and Laporte, G. Improved formulations and algorithmic components for the electric vehicle routing problem with nonlinear charging functions. *Computers and Operations Research*, 2019. 104:256–294.

24. Kancharla, S. R., and Ramadurai, G. Electric vehicle routing problem with non-linear charging and load-dependent discharging. *Expert Systems with Applications*, 2020. 160:113714.

25. Sweda, T. M., Dolinskaya, I. S., and Klabjan, D. Adaptive routing and recharging policies for electric vehicles. *Transportation Science*, 2017. 51(4):1326–1348.

26. Keskin, M., Laporte, G., and Çatay, B. Electric vehicle routing problem with time-dependent waiting times at recharging stations. *Computers and Operations Research*, 2019. 107:77–94.

27. Kullman, N. D., Goodson, J. C., and Mendoza, J. E. Electric vehicle routing with public charging stations. *Transportation Science*, 2021. 55(3):637–659.

28. Keskin, M., Çatay, B., and Laporte, G. A simulation-based heuristic for the electric vehicle routing problem with time windows and stochastic waiting times at recharging stations. *Computers and Operations Research*, 2021. 125:105060.

29. Keskin, M., and Çatay, B. Partial recharge strategies for the electric vehicle routing problem with time windows. *Transportation research part C: emerging technologies*, 2016. 65:111–127.

30. Schiffer, M., and Walther, G. The electric location routing problem with time windows and partial recharging. *European journal of operational research*, 2017. 260(3):995–1013.
22